\documentclass[preprint,11pt]{elsarticle}

\usepackage[latin1]{inputenc}
\usepackage[T1]{fontenc}
\usepackage{amssymb}
\usepackage{amsfonts}
\usepackage{amsmath}

\newtheorem{proposition}{Proposition}
\newtheorem{definition}{Definition}

\newcommand{\norm}[1]{\left\Vert{#1}\right\Vert}  

\newcommand{\sign}[1]{\mathrm{sign}\left(#1\right)}

\newcommand{\be}{\begin{equation}}
\newcommand{\ee}{\end{equation}}

\newcommand{\ps}[2]{\langle #1,#2\rangle}

\newtheorem{theorem}{Theorem}
\newtheorem{lemma}{Lemma}
\newtheorem{corollary}{Corollary}
\newcommand{\xs}{x^{\star}}
\newcommand{\R}{\mathbb{R}}
\newenvironment{preuve}{\noindent}{\hfill$\square$\bigbreak} 
\title{Consistency of $\ell_1$ recovery from noisy deterministic measurements.}
\date{November 2012}
\author{Charles~Dossal}
\ead{charles.dossal@math.u-bordeaux1.fr}
\address{IMB, Universit\'{e} Bordeaux 1, 351, Cours de la
  Lib\'{e}ration, 33405 Talence Cedex, France}
\author{Remi~Tesson}
\ead{remi.tesson@eleves.bretagne.ens-cachan.fr}
\address{ENS Cachan - Bretagne
Campus de Ker Lann
Avenue Robert Schuman
35170 Bruz}
\begin{document}
\begin{abstract}
In this paper a new result of recovery of sparse vectors from
deterministic and noisy measurements by $\ell_1$ minimization is
given. The sparse vector is randomly chosen and follows a 
{\it generic $p$-sparse model} introduced by Candes and
al. \cite{CandesPlan09}. 
The main theorem ensures consistency of $\ell_1$ minimization
with high probability. 
This first result is secondly extended to compressible vectors.    
\end{abstract}
\maketitle
\section*{Introduction}
Let $A$ be a real matrix with $n$ rows and $m$ columns with $m>n$. 
Let $x^0$ be a sparse vector following a {\it generic p-sparse model} and 
let $y$ be a data vector $y=Ax^0+b$, where $b$ is a noise vector. 
The question we want to address is: can we give a bound on the
sparsity 
of $x^0$ ensuring 
$x^0$ can be recovered or estimated from $\ell_1$ minimization 
with high probability ?

Candès and Plan \cite{CandesPlan09} answer partially to this 
question under assumptions on the coherence of the matrix $A$, with a
random (Gaussian) noise and with hypotheses on the minimum 
absolute value of non-zero components of $x^0$. They proved that 
with high probability the support and the sign of $x^0$ can 
be recovered using $\ell_1$ minimization if $x^0$ is sparse enough. \\
In this paper we show that, under the same assumption on the 
coherence and sparsity, with a bounded noise, without any assumption 
on the minimum absolute 
value of $x^0$, $\ell_1$ minimization provides a vector 
$x^\star$ such that $\norm{x^0-x^\star}_2$ can be bounded. 
Moreover this new result can be extended to compressible 
vectors that are close to sparse vectors.\\
Since an explicite formulation of $\xs$ 
is impossible without the minimum value assumption,
different tools must be developed.

In a first part, notations and definitions are given. 
In a second part the contributions of the paper are developed and
connected to prior works. 
In a third part the proof of main results are given. 
A last part is devoted to discussion.   
\section{Notations and Definitions}

Let us recall the definition of the subgradient of the $\ell_1$ norm at a point $x$ which support is $I$ :
\begin{equation*}
\partial \norm{x}_1=\{\xi\text{ such that } \norm{\xi}_{\infty}\leqslant 1,\,\forall i\in I, \xi(i)=\sign{x(i)}\}.
\end{equation*}
The Bregman distance is defined as follows:
\begin{definition}
Let $x^1,\,x$ be two vectors of $\R^m$. For all $\xi\in\partial \norm{x^1}_1$, the Bregman distance between $x$ and $x^1$ is defined by
\begin{equation*}
D_{\xi}(x,x^1))=\norm{x}_1-\norm{x^1}_1-\ps{\xi}{x-x^1}.
\end{equation*}  
\end{definition}
The {\it generic p-sparse model} is defined by Candès and Plan \cite{CandesPlan09} as follows
\begin{definition}
A vector $x$ follows the {\it generic p-sparse model} if the support of $x$ is randomly chosen with equiprobability from all supports which cardinal is
$p$ and if its sign is randomly chosen with equiprobability from all possible sign vectors .
\end{definition}
\begin{definition}
For a matrix $B$ the norm $\norm{B}_{p\to q}$ is defined as follows :
\begin{equation*}
\norm{B}_{p\to q}=\underset{x\neq 0}{\sup}\frac{\norm{Bx}_q}{\norm{x}_p}\quad \text{and}\quad \norm{B}_p=\norm{B}_{p\to p}.
\end{equation*}
\end{definition}
For a matrix $B$ such that $B^tB$ is invertible, $B^+=(B^tB)^{-1}B^t$ denotes the Moore-Penrose pseudoinverse of $B$.\\
Recall that $\norm{B^+}_2=\sqrt{\norm{(B^tB)^{-1}}_2}$ and that $\norm{B}_{1\to 2}=\underset{i}{\max}\norm{b_i}_2$ where $(b_i)_i$ are the columns 
of $B$. 
\begin{definition}
A $n\times m$ matrix $A$ whose columns are normalized is said to satisfy the $A_0-$coherence criterion if  
\begin{equation}
\mu(A)=\underset{i\neq j}{\max}|\ps{a_j}{a_i}|\leqslant \frac{A_0}{\ln{m}}
\end{equation} 
where $A_0$ is non negative real number.
\end{definition}

For a given vector $x^0\in\R^m$ which support is $I$, $A_I$ denotes the submatrix of $A$ which columns are columns of $A$ indexed by $I$. For a given vector $x$, $x_I$ denotes the subvector which components are components of $x$ indexed by $I$. The vector $\sign{x}$ is the vector whose component indexed by $i$
is 1 if $x(i)>0$, and $-1$ if $x(i)<0$ and  0 if $x(i)=0$.
If the columns $(a_i)_{i\in I}$ are linearly indepedent, the matrix $A_I^tA_I$ is invertible and for any $x^0$ such that $Supp(x^0)=I$, one can define
\begin{equation*}
 d(x^0)=A_I(A_I^tA_I)^{-1}\sign{x^0_I}\quad\text{ and }IC(x^0)=\underset{j\notin I}{\max}|a_j^td(x^0)|.
\end{equation*}   
This Identification Coefficient ($IC$) can be seen as signed ERC 
(Exact Recovery Coefficient
introduced by Tropp \cite{tropp-just-relax,tropp-greed-good}).
The condition $IC(x^0)<1$, see Fuchs \cite{Fuchs-redundant-bases}, is a
sufficient condition for exact recovery by $\ell_1$
minimization. 

\section{Contributions and relations with prior works}
Suppose $y=Ax^0+b$ with $\norm{b}_2\leqslant \varepsilon$ and define the minimization problem
\begin{equation}\label{Pb1}
\underset{x\in \R^m}{\min}\norm{x}_1\quad\text{under the constraint}\quad \norm{Ax-y}_2\leqslant \varepsilon.
\end{equation}   
Let $\xs$ be a minimizer of \eqref{Pb1}. 

The following Theorem holds. 
\begin{theorem}\label{Th1}
Let $A$ be a $n\times m$ matrix 
satisfying the $A_0-$coherence criterion with $A_0$ small enough, 
suppose that $x^0$ follows the {\it generic p-sparse model}, with $p\leqslant \frac{c_0m}{\norm{A}_2^2\ln{m}}$, for $c_0$ small enough depending on $A_0$.
Suppose $\norm{b}_2\leqslant \varepsilon$ and $y=Ax^0+b$, then any solution $\xs$ of \eqref{Pb1} satisfies
\begin{equation}
\norm{x^0-\xs}_2\leqslant C\varepsilon
\end{equation}  
with 
\begin{equation}\label{BorneBruit}
C=2\sqrt{2}+\frac{8(2+\sqrt{2})\sqrt{p}}{3}.
\end{equation}
with probability greater than $1-4m^{-2\ln 2}$ if $m$ is large enough.
\end{theorem}
It turns out that $\norm{\xs-x^0}_1$ can also be bounded using a similar proof :
\begin{equation}
\norm{\xs-x^0}_1\leqslant \frac{\varepsilon}{3}(14\sqrt{2p}+16p),
\end{equation}
with the same probability. The proof of this inequality requires a simple modification of Proposition \ref{Pnorme2}. This extension may be interesting since vectors $x^0$ and $\xs$ belong the a space which dimension is $m$ much larger than $p$.\\
Applying the 
theorem to $\xs=x^1$ and $b_1=Ar+b$ one obtains the following corollary : 
\begin{corollary}\label{Cor1}
Suppose $x^1=x^s+r$ where $\norm{r}_2\leqslant C_1\varepsilon_1$ 
and $x^s$ follows the {\it generic p-sparse model}, with 
$p\leqslant \frac{c_0m}{\norm{A}_2^2\ln{m}}$, for $c_0$ 
small enough and if $A$ satisfies the $A_0-$coherence 
criterion with $A_0$ small enough.
Suppose $\norm{b}_2\leqslant \varepsilon_1$ then any solution 
$\xs$ of \eqref{Pb1} with $y=Ax^1+b$ and 
$\varepsilon=\varepsilon_1(1+C_1\norm{A}_2)$ satisfies
\begin{equation}
\norm{x^s-\xs}_2\leqslant C\varepsilon 
\end{equation}  
with 
\begin{equation}
C=2\sqrt{2}+\frac{8(2+\sqrt{2})\sqrt{p}}{3}.
\end{equation}
with probability greater than $1-4m^{-2\ln 2}$ if $m$ is large enough.
\end{corollary}
To prove the corollary, one can apply the Theorem with $x^0=x^s$ and $b_1=Ar+b$.\\
This result sheds a new light on the understanding of the success of
$\ell_1$ minimization of the recovery of sparse and compressible
vectors from noisy deterministic measurements. No Restricted Isometry
Properties (RIP) \cite{candes-cras,candes-near-optimal} 
can be used here. The geometry of polytopes associated to $A$ 
(see Donoho \cite{donoho-polytopes})
seems hard to use
 and the classical bound derived by the coherence or the ERC 
\cite{tropp-just-relax} are too weak.
In \cite{CandesPlan09} authors propose an approach with a random model 
on the vector $x^0$. This work lies on concentration lemmas of
singular 
values of submatrices due to Tropp and on a explicit formulation of 
the solution of $\ell_1$ minimization, see Fuchs
\cite{Fuchs-bounded-noise}. 
This approach ensures the exact recovery of the support and the 
sign of the solution and needs conditions on the signal to noise 
ratio, decorrelation between the noise and the matrix and 
consequently can not be easily extended to compressible vector.\\
The present article focuses on the $\ell_2$ reconstruction error. 
In this new setting no signal to noise ratio, 
no independence between matrix and noise are needed and the 
result can be easily extended to compressible vectors. 
However this new approach doesn't give any informations on 
the support of the solution.
Unlike \cite{CandesPlan09}, the bound holds for any vector 
$x^1$ satisfying $\norm{y-Ax^1}_2\leqslant \varepsilon$ and 
$\norm{x^1}_1\leqslant \norm{x^0}_1$.\\
Following Grasmair et al. \cite{Grasmair11}, 
our approach uses Bregman distance to bound the part of the $\ell_1$
norm of $\xs$ that is not supported on the support $I$ of $x^0$. 
\section{Proof of Theorem \ref{Th1}.}
\subsection{Proof of Theorem \ref{Th1}.}
The proof lies on two properties, the first one bounds the $\ell_2$ error $x^0-\xs$ under the hypothesis that $IC(x^0)<1$ 
\begin{proposition}\label{Pnorme2}
Let $x^0\in\R^m$, whose support is $I$.
If $IC(x^0)<1$, then for any $\xs$ solution of \eqref{Pb1}, the following inequality holds
\begin{equation}\label{eqProp1}
\norm{\xs-x^0}_2\leqslant 2\varepsilon\left(\sqrt{\norm{(A_I^tA_I)^{-1}}_{2}}+\frac{\norm{d(x^0)}_2}{1-IC(x^0)}(\sqrt{\norm{(A_I^tA_I)^{-1}}_{2}}\norm{A_{I^c}}_{1\to 2}+1)\right)
\end{equation}
\end{proposition} 
The second proposition ensures that if $x^0$ follows the $p-$sparse model for $p$ small enough then with high probability $IC(x^0)<\frac{1}{4}$ and $\norm{(A_I^tA_I)^{-1}}_2\leqslant 2$ :
\begin{proposition}\label{PropBorneF}
Suppose $x^0$ follows the {\it generic p-sparse model} with $p\leqslant \frac{c_0m}{\norm{A}_2^2\ln m}$ for $c_0$ small enough
\begin{equation}
P\left(\left(\norm{(A_I^tA_I)^{-1}}_2\leqslant 2\right)\bigcap\left(IC(x^0)<t\right)\right)\geqslant 1-2m \exp\left(-\frac{t^2\ln m}{8c_0^2}\right)-2m^{-2\ln 2}.
\end{equation}
\end{proposition}

Choosing $c_0$ small enough in Proposition \ref{PropBorneF} yields 
\begin{equation}
P\left(\left(IC(x^0)\leqslant \frac{1}{4}\right)\bigcap\left(\norm{(A_I^tA_I)^{-1}}_2\leqslant 2\right)\right)\geqslant 1-4m^{-2\ln 2}.
\end{equation}
Moreover  
\begin{equation}
\norm{d(x^0)}_2^2=\ps{\sign{x^0_I}}{(A_I^tA_I)^{-1}\sign{x^0_I}}\leqslant \norm{(A_I^tA_I)^{-1}}_{2}p.
\end{equation}
It can be noticed that for any support $I$, if the columns of $A$ are normalized, 
$\norm{A_{I^c}}_{1\to 2}=1$.

Applying Proposition \ref{Pnorme2} to $x^0$, it follows that with probability greater than $1-4m^{-2\ln 2}$,
\begin{equation}
\norm{\xs-x^0}_2\leqslant \varepsilon \left(2\sqrt{2}+\frac{8(2+\sqrt{2})\sqrt{p}}{3}\right),
\end{equation}
which concludes the proof of Theorem \ref{Th1}.
\subsection{Proof of proposition \ref{Pnorme2}}
\begin{preuve} The proof of Proposition \ref{Pnorme2} follows the one of Grasmair
  et al. in \cite{Grasmair11} using the fact that, under the
  assumption $IC(x^0)<1$, $s=A^td(x^0)\in\partial \norm{x^0}_1$. 
\begin{align*}
\norm{\xs-x^0}_2&\leqslant \norm{\xs_I-x_I^0}_2+\norm{\xs_{I^c}}_2\\
&\leqslant \norm{A_I^+A_I(\xs_I-x^0_I)}_2+\norm{\xs_{I^c}}_1\\
& \leqslant \norm{A_I^+}_{2}\norm{A_I(\xs_I-x^0_I)}_2+\norm{\xs_{I^c}}_1\\
& \leqslant \norm{A_I^+}_{2}(2\varepsilon+\norm{A_{I^c}\xs_{I^c}}_2)+\norm{\xs_{I^c}}_1\\
& \leqslant \norm{A_I^+}_{2}(2\varepsilon+\norm{A_{I^c}}_{1\to 2}\norm{\xs_{I^c}}_1)+\norm{\xs_{I^c}}_1\\
& \leqslant
2\varepsilon\norm{A_I^+}_{2}+(\norm{A_I^+}_{2}\norm{A_{I^c}}_{1\to 2}+1)\norm{\xs_{I^c}}_1
\end{align*}
Using the Bregman distance, $\norm{\xs_{I^c}}_1$ can be bounded : indeed, 
from the definition of $s=A^td(x^0)$ it follows that 
\begin{align*}
D_s(\xs,x^0)&=\norm{\xs}_1-\norm{x^0}_1-\ps{s}{\xs-x^0}\\
&=\norm{\xs}_1-\ps{s}{\xs}\\
&=\sum_{i\in I}(\sign{\xs_i}-\sign{x^0_i})\xs_i+\sum_{j\notin I}(\sign{\xs_j}-s_j)\xs_j\\
&\geqslant \sum_{j\notin I}(\sign{\xs_j}-s_j)\xs_j.
\end{align*}
Since $\forall j\notin I,\,|s_j|\leqslant IC(x^0)$, one gets
\begin{equation*}
D_s(\xs,x^0)\geqslant \sum_{j\notin I}(1-IC(x^0))\sign{\xs_j}\xs_j=(1-IC(x^0))\norm{\xs_{I^c}}_1,
\end{equation*}
that is, 
\begin{equation*}
\norm{\xs_{I^c}}_1\leqslant \frac{D_s(\xs,x^0)}{1-IC(x^0)}.
\end{equation*}
Consequently,
\begin{equation}\label{eqProp11}
\norm{\xs-x^0}_2\leqslant 2\varepsilon\norm{A_I^+}_{2}+(\norm{A_I^+}_{2}\norm{A_{I^c}}_{1\to 2}+1)\frac{D_s(\xs,x^0)}{1-IC(x^0)}.
\end{equation}
The Bregman distance can be bounded as follows. Since $IC(x_0)<1$, $s=A^td(x^0)\in\partial \norm{x^0}_1$.
Consequently
\begin{align*}
D_s(\xs,x^0)&=\norm{\xs}_1-\norm{x^0}_1-\ps{A^td(x^0)}{\xs-x^0}\\
&\leqslant -\ps{A^td(x^0)}{\xs-x^0}\\
&=-\ps{d(x^0)}{A(\xs-x^0)}\\
&\leqslant \norm{d(x^0)}_2\norm{A(\xs-x^0)}_2\\
&\leqslant \norm{d(x^0)}_2(\norm{A\xs-y}_2+\norm{b}_2)\\
&\leqslant 2\norm{d(x^0)}_2\varepsilon.
\end{align*}
The fact that $\norm{A_I^+}_2=\sqrt{\norm{(A_I^tA_I)^{-1}}_2}$ concludes the proof of Proposition \ref{Pnorme2}. 
\end{preuve}
\subsection{Proof of Proposition \ref{PropBorneF}}
\begin{preuve}
The proof relies on a proposition due to Tropp
\cite{tropp-NormsRandom} (see also \cite{CandesPlan09}):
\begin{proposition}\label{PropTropp}
Suppose that the set $I$ is randomly and uniformely chosen among sets of cardinal $p$ with 
$p\leqslant \frac{m}{4\norm{A}_2^2}$. Then for $q=2\ln m$, 
\begin{equation}
E(\norm{A_I^tA_I-Id}_2^q)^{\frac{1}{q}}\leqslant 30\mu(A)\ln m+13\sqrt{\frac{2p\norm{A}_2^2\ln m}{m}}
\end{equation}
and 
\begin{equation}
E(\underset{j\notin I}{\max}\norm{A_I^ta_j}_2^q)^{\frac{1}{q}}\leqslant 4\mu(A)\sqrt{\ln m}+\sqrt{\frac{p\norm{A}_2^2}{m}}.
\end{equation}
\end{proposition}
From this proposition and Markov inequality, Candès and Plan \cite{CandesPlan09} proved the following Corollary :
\begin{corollary}\label{CorTropp}
Suppose that $A$ satisfies the $A_0$-coherence criterion and that $x^0$ follows the {\it generic p-sparse model} with $p\leqslant \frac{c_0m}{\norm{A}_2^2\ln m}$ and $30A_0+13\sqrt{2c_0}\leqslant \frac{1}{4}$.
Then $A_I^tA_I$ is invertible with probability greater than $1-m^{-2\ln 2}$  and 
\begin{equation}\label{eqvpai}
\norm{(A_I^tA_I)^{-1}}_2\leqslant 2.
\end{equation}
\end{corollary}
From the Proposition \ref{PropTropp} and the Hoeffding inequality the following Lemma can be deduced  (see Candès and Plan \cite{CandesPlan09}):
\begin{lemma}\label{LW}
Suppose $x$ follows the {\it generic p-sparse model} and let $(W_j)_{j\in J}$ be a collection of deterministic vectors. For $Z_0=\underset{j\in J}{\max}|\ps{W_j}{\sign{x_I}}|$ one has
\begin{equation*}
P(Z_0\geqslant t)\leqslant 2|J|e^{-\frac{t^2}{2\kappa^2}}
\end{equation*} 
for $\kappa\geqslant \underset{j\in J}{\max}\norm{W_j}_2$.
\end{lemma}
Applying Lemma \ref{LW} and Proposition \ref{PropTropp}, the proof of Proposition \ref{Pnorme2} can be achieved :

For all $j\notin I$, define $W_j=a_j^tA_I(A_I^tA_I)^{-1}$, where $I$ is the support of $x^0$. Applying
Lemma \ref{LW}, one gets 
\begin{equation}\label{eqBorneFuchs}
P(IC(x)\geqslant t)\leqslant 2|I^c|e^{-\frac{t^2}{2\kappa^2}}\leqslant 2me^{-\frac{t^2}{2\kappa^2}}.
\end{equation}
We need to estimate the maximum of $\norm{W_j}_2$. Using Corollary \ref{CorTropp} one gets
\begin{equation}\label{eqW1}
\norm{W_j}_2=\norm{(A_I^tA_I)^{-1}A_I^ta_j}_2\leqslant \norm{(A_I^tA_I)^{-1}}_2\norm{A_I^ta_j}_2\leqslant 2\norm{A_I^ta_j}_2
\end{equation}
with a probability greater than $1-m^{-2\ln 2}$. Proposition \ref{PropTropp} and Markov inequality is then used to estimate $\norm{A_I^ta_j}_2$.
\begin{align*}
P\left(\underset{j\notin I}{\max}\norm{A_I^ta_j}_2>\frac{c_0}{\sqrt{\ln m}}\right)
&\leqslant c_0^q\frac{E(\underset{j\notin I}{\max}\norm{A_I^ta_j}_2^q)}{(\ln m)^{\frac{q}{2}}}\\
&\leqslant \left(\frac{c_0}{\sqrt{\ln m}}\right)^q\left(4\mu(A)\sqrt{\ln m}+\sqrt{\frac{p\norm{A}_2^2}{m}}\right)^q\\
&\leqslant \left(\frac{c_0}{\sqrt{\ln m}}\right)^q\left(\frac{4A_0+\sqrt{c_0}}{\sqrt{\ln m}}\right)^q.
\end{align*}
If $A_0$ and $c_0$ are small enough, $\displaystyle \frac{4A_0+\sqrt{c_0}}{\sqrt{\ln m}}\leqslant \frac{c_0}{2\sqrt{\ln m}}$ and 
using $q=2\ln m$ it follows
\begin{equation}\label{eqW2}
P\left(\underset{j\notin I}{\max}\norm{A_I^ta_j}_2>\frac{c_0}{\sqrt{\ln m}}\right)\leqslant m^{-2\ln 2}.
\end{equation}
From \eqref{eqW1} and \eqref{eqW2} it follows
\begin{equation}
P\left(\underset{j\notin I}{\max}\norm{W_j}_2>\frac{2c_0}{\sqrt{\ln m}}\right)\leqslant 2m^{-2\ln 2}.
\end{equation}
Combined with inequality \eqref{eqBorneFuchs}, this last inequality concludes the proof of the proposition.
\end{preuve}
\section{Discussions}
The two constants $\frac{1}{4}$ and $2$ in the proof of Theorem \ref{Th1} are arbitrary chosen. 
If others bounds are chosen, the optimal value of $c_0$ changes and the 
value of $C$ in Theorem \ref{Th1} may change. It turns out that these 
values are numerically pessimistic and that their optimization would not be useful. 
Two relevant questions may be asked about Theorem \ref{Th1} : 

Can we expect better bounds 
on the sparsity using the criterion $IC<1$? 

Constants may be optimized but 
it seems that the asymptotic of the sparsity may not be improved.
 In \cite{2012-dossal-acha} and \cite{wainwright-sharp-thresh} authors 
 proved that for gaussian measurements, beyond sparsity $\frac{m}{2\ln m}$, with high probability,
 $IC(x^0)>1$. It would be surprising that better results could be achieved by deterministic measurements.
    
 However the Grasmair approach (see Proposition \ref{Pnorme2}) applies 
 also to any vector $\eta$ in the subgradient of the $\ell_1$ norm at
 the point $x^0$ not only to $s=A^td(x^0)$. It may possible to improve
 the sparisty bound using another vector $\eta$. 
 
 The second question is about the $\sqrt{p}$ scaling in the bound \eqref{BorneBruit}.
 Is this scaling optimal or not ? Can we expect a better bound ?
 
 RIP Theory gives similar bounds where the constant $C$ in \eqref{BorneBruit} does not depend on the sparsity $p$ but on RIP constants that can be uniformely bounded if the vector is sparse enough. Moreover Fuchs \cite{Fuchs-bounded-noise} proved that when the noise is small enough and if $IC(x^0)<1$ the support of $\xs$ is equal to the support of $x^0$, that is $\norm{\xs_{I^c}}_1=0$. Looking at the proof of Proposition \ref{Pnorme2}, it appears that if $\xs_{I^c}=0$, in Theorem 
 \ref{Th1} the constant $C$ can be set to $2\sqrt{2}$ which do not depend on $p$. 
 Unfortunatly if no assumptions are made on $\varepsilon$, there is no guarantee that $x_{I^c}=0$. RIP theory solves the problem ensuring that all submatrices with a small number of columns have a good behaviour. Such hypothesis can not be done here and for some noise vectors $b$ 
 it may happen that $\norm{\xs_{I^c}}_1\neq 0$. If $I^{\star}$ denotes the support of $\xs$, the solution $\xs$ satisfies the following implicit equation (see \cite{Fuchs-bounded-noise})
 \begin{equation} \xs_{I^{\star}}=(A_{I^{\star}}^tA_{I^{\star}})^{-1}A_{I^{\star}}^ty-\lambda(A_{I^{\star}}^tA_{I^{\star}})^{-1}\sign{x_{I^{\star}}}
 \end{equation}  
 where $\lambda$ depends on $\varepsilon$. This expression shows that
 the stability of the solution depends widely on the matrix
 $(A_{I^{\star}}^tA_{I^{\star}})^{-1}$ wich depends on $x^0$ and on
 the noise $b$. In practice in many situations $I^{\star}\not\subset
 I$ and there is no simple way to control  $(A_{I^{\star}}^tA_{I^{\star}})^{-1}$.\\

The scaling $\sqrt{p}$ may be the price to pay of the lack of control on this matrix. 
  
\section{Conclusion}
These results complete the previous one of Candes and Plan
\cite{CandesPlan09} and ensures that under the same hypothesis of
sparsity $\ell_1$ minimization is robust to noise and compressibility
even if the exact support and sign can not be recovered. To controle
the part of the solution that is not supported on the support $I$ of the
objective vector $x^0$, no RIP can be used here but the Bregman
distance provides a interesting bound.    
\bibliographystyle{elsarticle-num}

\end{document}